\documentclass[12pt,reqno]{amsart}
\usepackage{pstricks}
\usepackage{amssymb, a4, mathrsfs}
\usepackage{amsmath}
\newpsobject{showgrid}{psgrid}{subgriddiv=1,griddots=5,gridlabels=6pt}
\usepackage[top=2.25cm, bottom=2.5cm, left=2.5cm, right=2.5cm ]{geometry}
\usepackage{paralist}
\usepackage{bbding}
\usepackage{graphicx}
\usepackage{url}
\usepackage{fancyhdr}
\usepackage{footmisc}
\usepackage{algorithmic}
\usepackage{listings}
\usepackage{fancyvrb}
\usepackage{cite}

\usepackage{tikz}

\usetikzlibrary{arrows.meta}%»­¼ýÍ·ÓÃµÄ°ü
\usepackage{amscd}

\usepackage{indentfirst}
\usetikzlibrary{arrows,automata}

\newtheorem{them}{Theorem}[section]
\newtheorem{prop}[them]{\noindent Proposition}
\newtheorem{lem}[them]{\noindent Lemma}
\newtheorem{coro}[them]{\noindent Corollary}

\theoremstyle{definition}
\newtheorem{definition}[them]{\noindent Definition}
\newtheorem{ex}[them]{\noindent Example}

\begin{document}
%\linenumbers  % ¿ªÆôÐÐºÅ

\title[Zero-divisor graphs of inverse semigroups]{On the diameter and girth of zero-divisor graphs of inverse semigroups}
\thanks{ }
\author{Yanhui Wang}
\address{College of Mathematics and Systems Science\\ Shandong University of Science and Technology\\ Qingdao 266590\\ China}
\email{yanhuiwang@sdust.edu.cn}
\author{Xinyi Zhu}
\address{College of Mathematics and Systems Science\\ Shandong University of Science and Technology\\ Qingdao 266590\\ China}
\email{zhuxinyi12333@163.com}
\author{Pei Gao}
\address{College of Mathematics and Systems Science\\ Shandong University of and Technology\\ Qingdao 266590\\ China}
\email{gaopei0403@163.com}

\noindent
\noindent
\maketitle

\renewcommand{\thefootnote}{\empty}

\vspace{3mm}
\begin{abstract}
Let $S$ be an inverse semigroup with zero and let $Z(S)^\times$ be its set of non-zero divisors with respect to the natural partial order $\le $ on $S$, that is, $a \in Z(S)^\times $ if there exists $b\in S\setminus\{0\}$ with $\omega(a, b) = \{c \in S: c \leq a\ \mbox{and}\ c \leq b\}=\{0\}$. The set $Z(S)^\times$ makes up the vertices of the corresponding {\it zero-divisor graph} $\Gamma (S)$, with two distinct vertices $a, b$ forming an edge if $\omega(a, b)=\{0\}$. We characterize {\it zero-divisor graphs} of inverse semigroups in terms of their diameter and girth. We also classify inverse semigroups without zero by  building a connection between the diameter (girth) and the least group congruence $\sigma$ on an inverse semigroup without zero. Finally, we give a description of the diameter and girth of graph inverse semigoups $I(G)$ in terms of the set of vertices and the set of edges of a graph $G$.
\end{abstract}

\section{Introduction}

Associating a graph to an algebraic structure is a research subject in algebraic combinatorics  and has attracted considerable attention~\cite{ref19, ref12,Chattopadhyay,Das,Jafaria,ref14}. Zero-divisor graphs play an important role in exposing the relationship between algebra and graph theory~\cite{ref19,ref12,Sbarra}.  There are primarily two ways to define the zero-divisor graph: one is based on the operations of algebraic systems, and the other is based on the order structure. Beck~\cite{ref5}  introduced the notion of a zero-divisor graph $\Gamma_0(R)$ of a commutative ring $R$ with identity to be the undirected graph whose vertices are elements of $R$ and in which two vertices $x$ and $y$ are adjacent if and only if $xy = 0$, where $xy$ is the product of $x$ and $y$ in $R$.   Many authors have studied zero-divisor graphs of rings~\cite{ref7} or the other
algebraic structures such as posets\cite{ref9} to show that Beck's conjecture, that is, $\chi(R) = \omega(R)$, where $\chi(R)$ and $\omega(R)$ denote
the chromatic number and the clique number of the zero-divisor groups $\Gamma_0(R)$, respectively.

 More recently, a different method of associating a zero-divisor graph to a poset $(P, \leq)$ was proposed by Lu and Wu in~\cite{ref11} using the partial order $\leq$.  Let $(P, \leq)$ be a poset with a least element $0$ and $P^{\times} = P\setminus\{0\}$. For all $x, y \in P$, $\omega(x, y) = \{z \in P: z \leq x\ \mbox{and}\ z \leq y\}$. The {\it zero-divisor graph} $\Gamma$  of a poset $P$ is an undirected graph  consists of a set $V$ of vertices and a set $E$ of edges, where $V = \{x \in P^{\times}:  \omega(x, y)=\{0\} \ \mbox{for\ some}\ y \in P^{\times} \}$ and for all $x, y \in V$, $x$ and $y$ are {\it adjacent} in $\Gamma$, that is, $\{x, y\} \in E$, if $\omega(x, y) = \{0\}$. Alizadeh et al. in~\cite{ref13} proved that the diameter of the zero-divisor graph associated with a poset is either $1$, $2$ or $3$ while its girth is either $3$, $4$ or $\infty$, and also classified
zero-divisor graphs of posets in terms of their diameter and girth.

 Mitsch showed that for an arbitrary semigroup $S$ there exists a partial order $\leq$ associated with it,  where $\leq$ is defined by means of the multiplication of $S$~\cite{Mitsch}.  An interesting question is that:

{\it If we deal with zero-divisor graphs of semigroups based on the terminology of ~\cite{ref11}, how to characterize zero-divisor graphs of semigroups in terms of their diameter and girth?}

In this paper, we give an answer to the above question for inverse semigroups.  The structure of this paper is organized as follows. Section \ref{sec:Pre} recalls some basic definitions and notation related to graphs and  inverse semigroups. Section~\ref{sec:zero-divisor} gives the definition of the {\it zero-divisor graph} of an inverse semigroup and classified the {\it zero-divisor graphs} of inverse semigroups in terms of their diameter and girth as an application of ~\cite[Theorem 3.3]{ref13} and ~\cite[Theorem 4.2]{ref13}.  Section \ref{proper semigroups}  considers  the  necessary and sufficient  conditions for any two vertices to be connected in zero-divisor graphs associated with inverse semigroups adding an extra  zero and also we give necessary and sufficient conditions when the diameter and the girth of its zero-divisor graph taking a value. Theorem~\ref{GammaS0} shows that the zero-divisor graph of an inverse semigroup $S^0$ with zero is close related to the least group congruence  $\sigma$ on $S$, that is, for $a, b \in S$, $a\, \sigma\, b$ if and only if $ea = fb$ for some $e, f \in E$. In Section~\ref{graph inverse semigroup}, we are interested in   zero-divisor graphs of graph inverse semigroups. Theorem \ref{complete} indicates that $\mbox{Path}(G)\cup\mbox{Path}(G)^{-1}$ induced a complete subgraph of $\Gamma(I(G))$, and Proposition \ref{vertex set} shows that $V(\Gamma(I(G)))=I(G\setminus\{0\})$. Finally we show that the diameter of the zero-divisor graph associated with a graph inverse semigroup is either $1$ or $2$ while its girth is either $3$ or $\infty$.

\section{Preliminaries}\label{sec:Pre}
To make this article self-contained we recall some basic definitions and properties concerning graphs and inverse semigroups. For more details, we refer the reader to \cite{ref16} and \cite{ref1}.
		
\subsection{Graphs}
An {\it undirected graph} $G = (V(G), E(G))$ consists of a set $V(G)$ of vertices and a set $E(G)$ of edges.  Any two vertices $u$ and $v$  in $G$ are {\it adjacent} if there exists an edge $e \in E(G)$ such that $u$ and $v$ are two endpoints of $e$. A directed graph $G = (V(G), E(G), {\bf s}, {\bf r})$ consists of a {\it set of vertices} $V(G)$, a {\it set of edges} $E(G)$, and two mappings $\mathbf{s}, \mathbf{r} : E(G) \rightarrow V(G)$, respectively, called the {\it source mapping} and the {\it range mapping} for $G$. If an edge $e$ is starting at $u$ and ending at $v$, then  we write as ${\bf{s}}(e) = u$  and ${\bf{r}}(e) = v$, respectively. A graph $G$ with only a few isolated vertices is called a {\it null graph} and specially $G$ is a {\it trivial graph} if it has one isolated vertex. Throughout this paper we will explicitly mention when we consider directed graphs, otherwise ``graph''  will refer to a {\it simple} undirected graph that is an undirected graph without loops and multiple edges.

A {\it path} $p$ of length $r$ from $u$ to $v$ in a graph  is a sequence of $r+1$ vertices starting at $u$ and ending at $v$ such that consecutive vertices are adjacent. Here $u$ is called the {\it source} of $p$ and $v$ is called the {\it range} of $p$ in a directed graph. We write as ${\bf{s}}(p) = u$ and  ${\bf{r}}(p) = v$. We also consider a vertex $v$ as being an {\it empty} path (i.e. a path with no edges) based at $v$ and with
${\bf{s}}(v)= {\bf{r}}(v) = v$. We denote the set of paths of $G$ by Path$(G)$.

It is convenient to extend the notation so as to allow
paths in which edges are read in either the directed or
inverse direction. To do this, we associate with each edge $e$ an
``inverse edge'' $e^{-1}$ (sometimes called a ``ghost edge'' by some
authors) with $\mathbf{s}(e^{-1}) = \mathbf{r}(e)$ and $\mathbf{r}(e^{-1}) = \mathbf{s}(e)$. We denote by $E(G)^{-1}$ the set $\{e^{-1} : e \in E(G)\}$ and assume that $E(G) \cap E(G)^{-1} = \emptyset$. With this convention, for each path $p =
e_1e_2\ldots e_n$ in $G$ we have an {\it inverse path} $p^{-1}
= e_n^{-1} \ldots e_2^{-1} e_1^{-1}$ and vice versa. As usual, $\mathbf{s}(p^{-1}) = \mathbf{s}(e_n^{-1}) = \mathbf{r}(e_n)$ and $\mathbf{r}(p) = \mathbf{r}(e_1^{-1}) = \mathbf{s}(e_1)$. We put $\mbox{Path}(G)^{-1}=\{p^{-1}: p\in \mbox{Path}(G)\}$.

It is easy to see that the {\it length} of a path is the number of edges in the path.  A {\it cycle} of a graph is a path such that the start and end vertices are the same. We refer to a cycle with $k$ edges as a {\it $k$-cycle}. If $k$ is odd, we call a k\mbox{-}cycle  an {\it odd cycle}. The {\it distance} between $u$ and $v$ in a graph $G$, denoted by ${\bf d}(u, v)$, is the length of a shortest path connecting $u$ and $v$, where $u$ and $v$ are  distinct vertices in $G$. If there is no any path between $u$ and $v$, we write ${\bf d}(u, v)=\infty$. The largest distance   among all distances between pairs of vertices of a graph $G$ is called the {\it diameter} of $G$ and is denoted by ${ \mbox{\bf diam}}(G)$. The {\it girth} of $G$ is the length of a shortest cycle in $G$ and is denoted by $\mbox{\bf gr} (G)$. If $G$ has no cycles, we define the girth of $G$ to be infinite.

A {\it subgraph} of a graph $G$ is a graph $G'$ such that $V(G') \subseteq V(G)$ and $E(G') \subseteq E(G)$. A subgraph $G'$ of $G$ is an {\it induced subgraph} if two vertices of $V(G')$ are adjacent in $G'$ if and only if they are adjacent in $G$. If there is a path between any two vertices of a graph $G$, then $G$ is {\it connected}, otherwise {\it disconnected}.

 A {\it bipartite graph} is one whose vertex-set is partitioned into two disjoint subsets in such a way that the two endpoints for each edge lie in distinct partitions. Among bipartite graphs, a {\it complete bipartite graph} is one in which each vertex is joined to every vertex that is not in the same partition. The complete bipartite graph with exactly two partitions of size $m$ and $n$ is denoted by $K_{m,n}$. Graphs of the
form $K_{1,n}$ are called {\it star graphs}. A graph is called {\it complete} if every pair of vertices are adjacent.

\subsection{Inverse semigroups}
In this subsection, we recall some results and properties about  inverse semigroups.
	
Let $S$ be a semigroup. An element $e \in S$ is called an {\it idempotent} if $e^2=e$. We denote the set of all idempotents of $S$ by $E(S)$.  If $S$ with at least two elements contains an element $0$ such that, for all $x$ in $S$,
$$0x=x0=0,$$
we say that $0$ is a {\it zero element} (or just a zero) of $S$, and that $S$ is a {\it semigroup with zero}. It is easy to see that there can be at most one zero element in a semigroup. And $0$ is also an idempotent of $S$. We use $S^0$ to denote $S$ with an external  zero element 0 adjacent if $0 \notin S$, otherwise $S^0=S$. If $S$ is a semigroup with a zero element $0$, we denote $S\setminus\{0\}$ by $S^{\times}$. In particular, if $0 \notin S$ we have $S^{\times}=S$.

A semigroup $S$ is said to be an {\it inverse semigroup} if, for every $a$ in $S$, there exists an element $b$ in $S$ such that $a = aba$ and $b = bab$, and $E(S)$ is a semilattice, where $E(s)$ is the set of idempotents of $S$.  V. Vagner defined a natural partial order on an inverse semigroup $S$ as follows:
\[a \leq b \ \mbox{if\ and\ only\ if}\ a = eb\ \mbox{for\ some}\ e \in E(S).\]
which is equivalent to the following  statement
\[a \leq b \ \mbox{if\ and\ only\ if}\ a = bf\ \mbox{for\ some}\ f \in E(S).\]

Let $S$ be an inverse semigroup with semilattice $E(S)$. For all $a \in S$, we set
$$\omega (a)=\{x \in S: x\le a\}.$$
 Then it is clear that
 $$\omega (a)=\{x \in S: x\le a\}=Ea=aE.$$

 An element  $b \in S$ is called a {\it minimal} element of $S$ if $x \in S$ and $x \leq b$ implies that $x = b$. We denote the set of minimal elements of $S$ by $\mbox{Min}(S)$.  For every $x \in S$, if there exists $a \in S$ such that $a \leq x$ , then $a$ is called   the {\it least} element of $S$.  The least element, if exists, is unique because of the antisymmetry of the partial order. If $S$ contains a zero $0$, for all $x \in {S}$, we have $0=0x$ namely $0 \leq {x}$ which implies that $0$ is the least element of $S$. The converse is true as follows:

\begin{lem}\label{0 inverse}
Let $S$ be an inverse semigroup with semilattice $E(S)$. Then an element $a \in S$ is the least element of $S$ with respect to the natural partial order if and only if $a$ is the zero element $0$ of $S$.
\end{lem}
\begin{proof}
It is sufficient to show that if $a$ is the least element of $S$ then $a$ is the zero element. Suppose that $a$ is the least element of $S$. Then for all $e \in E(S)$, we have $a \leq e$, that is, $a = fe$ for some $f \in E(S)$,  then $a \in E(S)$ and $ae = ea = a$ as $E(S)$ is a semilatttice. For all $b \in S$, we have $a \leq b$ and then $a = gb=bh$ for some $g,h \in E(S)$. Further, we have $ab = (ag)b = a(gb)=a^2 = a$ and $ba =b(ha)=(bh)a=a^2 = a$. Hence $a$ is a zero element of $S$.
\end{proof}

In this paper we consider the partial order on an inverse semigroup is the natural partial order.

\section{The zero-divisor graph}\label{sec:zero-divisor}

In this section, we define the {\it zero-divisor graph} of an inverse semigroup with zero and characterize its diameter and girth.

Let $S$ be an inverse semigroup with a zero element $0$. For every $a, b\in S$, denote
$$\omega(a, b)=\{c\in S: c\leq a\ \mbox{and}\ c\leq b\}.$$
 For $x \in S$, the {\it annihilator} of $x$, denoted by Ann$(x)$, is defined to be $\{y \in S: \omega(x, y) = \{0\}\}$.

\begin{definition}\label{ZeroDivisor}
 An element $a \in S$ is called a {\it zero-divisor} of $S$ if there exists $b \in S^\times $ such that $\omega(a, b)=
     \{0\}$.
 \end{definition}

 Set
 $$Z(S)=\{a\in S: \exists b \in S^\times, \omega(a, b)=
     \{0\}\}.$$
     and  $Z(S)^{\times} = Z(S)\setminus\{0\}$.

\begin{definition}\label{zero-divisor}
Let $S$ be an inverse semigroup with a zero element $0$.  The {\it zero-divisor graph} of $S$, denoted by $\Gamma (S)$, is the graph  whose set of vertices is  $Z(S)^{\times}$ and two distinct vertices $a$ and $b$ are adjacent if $\omega(a, b)=\{0\}$.
\end{definition}

Example~\ref{ex1} indicates that  the zero-divisor graph  defined in Definition~\ref{zero-divisor} are distinct from the one in a ring $R$. An element $a$ of a ring $R$ is called a {\it left zero-divisor} if there exists a nonzero $x$ in $R$ such that $ax = 0$~\cite{ref18}. Similarly, an element $a$ of a ring is called a {\it right zero-divisor} if there exists a nonzero $y$ in $R$ such that $ya = 0$. An element that is a left or a right zero divisor is simply called a {\it zero divisor}~\cite{ref17}.

\begin{ex}\label{ex1}	Let $S=B_2=\{a,b:a^2=b^2=0,\ aba=a,\ bab=b\}$. We denote $ab$ and $ba$ by $e$ and $f$, respectively. The Caylay table of $S$ is as follows.
	$$
	\begin{array}{l|lllll}
		* & 0 & e & f & a & b     \\ \hline
		0 & 0 & 0 & 0 & 0 & 0 \\
		e & 0 & e & 0 & a & 0    \\
		f & 0 & 0 & f & 0 & b   \\
		a & 0 & 0 & a & 0 & e    \\
		b & 0 & b & 0 & f & 0 \ . \\
	\end{array}
	$$
It is easy to see that $E(S)=\{0, e, f\}$ and $S$ is an inverse semigroup. Notice that $\omega(a, b)=\{0\}$, so $a$ and $b$ are adjacent in $\Gamma(S)$. Since $a^2=0$ and $b^2=0$ it follows that $a$ and $b$ are zero-divisors according to the definition of zero-divisors in a ring. But $ab=e\ne 0$ and $ba=f\ne 0$, then $a$ and $b$ are not adjacent in the zero-divisor graph defined as that of a ring\cite{ref19}.
\end{ex}

It is easy to obtain the following two propositions so we omit its proof.
\begin{prop}\label{a}
Let $S$ be an inverse semigroup with a zero element $0$. For all $a, b \in S^{\times}$, the following statements are equivalent:
\begin{enumerate}
\item[\rm (i)] $a$ and $b$ are adjacent in $\Gamma(S)$;
\item[\rm (ii)] $\omega(a, b) = \{0\}$;
\item[\rm (iii)] $Ea\cap Eb=\{0\}$;
\item[\rm (iv)]$aE\cap bE=\{0\}$;
\item[\rm (v)] $a \in \text{Ann}(b)$;
\item[\rm (vi)] $b \in \text{Ann}(a)$.
\end{enumerate}
\end{prop}

\begin{prop} Let $S$ be an inverse semigroup with zero. We have \label{minimal}
 {\rm Min}$(S^{\times})={\mathcal{M}}$, where
 \[{\mathcal{M}}=\{x\in S\setminus\{0\}:Ex=\{x,0\}\}.\]
\end{prop}

A nontrivial inverse semigroup $S$ with zero together with the natural  partial order forms a poset.
By \cite{ref13}, we obtain the diameter and girth of the zero-divisor graph associated with poset (see \cite[Theorem 3.3]{ref13} and \cite[Theorem 4.2]{ref13}).  So as an application of \cite[Theorem 3.3]{ref13} and  \cite[Theorem 4.2]{ref13} we get the diameter and  girth of $\Gamma(S)$ as follows:

\begin{prop}\label{diamP}
Let $S$ be an inverse semigroup with zero. Then the following statements hold:
\begin{enumerate}
\item[\rm (i)] $\Gamma(S)$ is a connected graph with {\bf diam}$(\Gamma(S)) \in \{1, 2, 3\}$.
\item[\rm (ii)]  {\bf diam}$(\Gamma(S)) = 1$ if and only if $Z(S)^{\times} = \text{Min}(S^{\times})$.
\item[\rm (iii)]  {\bf diam}$(\Gamma(S)) = 2$ if and only if $Z(S)^{\times}\setminus \text{Min}(S^{\times}) \neq \emptyset$ and Ann$(x) \cap$ Ann$(y) \neq \{0\}$ for all $x, y \in Z(S)^{\times}\setminus \text{Min}(S^{\times})$ with $\omega(x, y) \neq \{0\}$.
\item[\rm (iv)] {\bf diam}$(\Gamma(S)) = 3$ if and only if $Z(S)^{\times}\setminus \text{Min}(S^{\times}) \neq \emptyset$ and Ann$(x) \cap$ Ann$(y) = \{0\}$ for some $x, y \in Z(S)^{\times}\setminus \text{Min}(S^{\times})$ with $\omega(x, y) \neq \{0\}$.
\end{enumerate}
\end{prop}

 \begin{prop}\label{girtt}
 	Let $S$ be an inverse semigroup with zero. Then the following statements hold:
 \begin{enumerate}
 	\item[\rm (i)]  $\mbox{gr}(\Gamma(S))\in \{3,4,\infty\}$.
 	\item[\rm (ii)]  $\mbox{gr}(\Gamma(S)) = \infty$ if and only if $\Gamma(S)$ is a star graph.
 	\item[\rm (iii)]  $\mbox{gr}(\Gamma(S)) = 4 $ if and only if $\Gamma(S)$ is a bipartite graph but not a star graph.
 	\item[\rm (iv)] $\mbox{gr}(\Gamma(S)) = 3 $ if and only if $\Gamma(S)$  contains an odd cycle.
 \end{enumerate}
 \end{prop}

\section{Inverse semigroups without zero}\label{proper semigroups}
The aim of this section is to describe the diameter and girth of {\it zero-divisor graph} of inverse semigroups using the least group congruence.

Let $S$ be an  inverse semigroup with semilattice of idempotents $E$. The relation $\sigma$ on $S$ defined by the rule that for all $a, b \in S$,
\[a\ \sigma\ b\ \mbox{if\ and\ only\ if}\ ea=fb\ \mbox{for\ some}\ e, f \in E,\]
is the least group congruence on $S$.

If $S$ is an inverse semigroup  with zero, we have $0 \in E$ and so for all $a, b \in S$, $0a = 0b =0$. It follows that $\sigma$ is the universal relation on $S$. So, in this case we can not use $\sigma$ to characterize the properties of {\it zero-divisor graphs} of inverse semigroups with zero. In the following we only consider inverse semigroups without zero.

\begin{lem}\label{properInv}
Let $S$ be an inverse  semigroup without zero and let $E$ be the semilattice of idempotent of $S$. For all $a, b \in S$, the following statements are equivalent:
\begin{enumerate}
\item[\rm (i)]  $(a, b) \notin \sigma$;
\item[\rm (ii)] $a$ and $b$ are  adjacent in $\Gamma(S^0)$;
\item[\rm (iii)]  $E^0a \cap E^0b = \{0\}$.
 \end{enumerate}
\end{lem}
\begin{proof}
 (i) $\Rightarrow$ (ii). Let $a, b \in S$ be such that $(a, b) \notin \sigma$. Suppose that $a$ and $b$ are not adjacent in $\Gamma(S^0)$, that is, $\omega(a, b) \neq \{0\}$. Then there exists $x \in S$ such that $x \leq a$ and $x \leq b$, that is, $x = ea = fb$ for some $e, f \in E$. It follows that  $a\ \sigma\ b$, a contradiction. Hence, we have $\omega(a, b)= \{0\}$, that is, $a$ and $b$ are  adjacent in $\Gamma(S^0)$.

According to Proposition~\ref{a}, we have (ii) $\Rightarrow$ (iii) and so we next show that (iii) $\Rightarrow$ (i). Suppose that $a, b \in S$ are such that $E^0a \cap E^0b = \{0\}$. Since $0 \notin S$ it follows that there does not exist $e, f \in E$ such that $ea = fb$, and so $(a, b) \notin \sigma$.
\end{proof}

According to Lemma~\ref{properInv} for $a,b \in S$ if $a\ne b$ and $a\ \sigma\ b$ then $a$ and $b$ are not adjacent in $\Gamma(S^0)$. So if there exists only one $\sigma$-class in an inverse semigroup then the zero-divisor graph is empty; the converse  is also true.  Hence  Proposition~\ref{VertexOfProIvs} is obtained directly.

\begin{prop}\label{VertexOfProIvs}
Let $S$ be an inverse semigroup without zero. Then we have
\begin{enumerate}
\item[\rm (i)] $|S/\sigma|=1$ if and only if $V(\Gamma(S^0)) = \emptyset$.
\item[\rm (ii)] $|S/\sigma| \geq 2$ if and only if $V(\Gamma(S^0)) = S$.
 \end{enumerate}
\end{prop}

Let $S$ be an inverse semigroup without zero. Suppose that $a, b \in S$ are such that $(a, b) \notin \sigma$. Put $U = a\sigma \cup b\sigma$. Clearly, $U$ is a disjoint union of $a\sigma$ and $b\sigma$. It follows from Lemma~\ref{properInv} that the induced subgraph $\Gamma(U)$ with set $U$ of vertices  is a complete bipartite graph two partitions of which are $\sigma$-classes: $a\sigma$ and $b\sigma$ as follows.

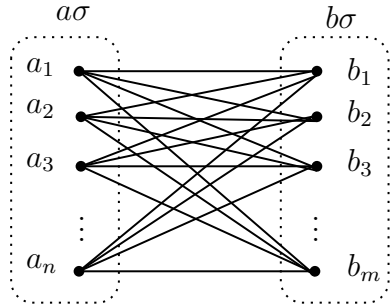
\tikzset{every picture/.style={line width=0.75pt}} %set default line width to 0.75pt
\begin{center}
\begin{tikzpicture}[x=0.75pt,y=0.75pt,yscale=-1,xscale=1]

%Rounded Rect [id:dp8977625277055314]
\draw  [dash pattern={on 0.84pt off 2.51pt}] (100,122.87) .. controls (100,116.87) and (104.87,112) .. (110.87,112) -- (143.47,112) .. controls (149.47,112) and (154.33,116.87) .. (154.33,122.87) -- (154.33,235.13) .. controls (154.33,241.13) and (149.47,246) .. (143.47,246) -- (110.87,246) .. controls (104.87,246) and (100,241.13) .. (100,235.13) -- cycle ;
%Shape: Circle [id:dp9450665838402348]
\draw  [fill={rgb, 255:red, 0; green, 0; blue, 0 }  ,fill opacity=1 ] (132,129.17) .. controls (132,127.95) and (132.99,126.96) .. (134.21,126.96) .. controls (135.43,126.96) and (136.42,127.95) .. (136.42,129.17) .. controls (136.42,130.39) and (135.43,131.37) .. (134.21,131.37) .. controls (132.99,131.37) and (132,130.39) .. (132,129.17) -- cycle ;
%Shape: Circle [id:dp8675857459590024]
\draw  [fill={rgb, 255:red, 0; green, 0; blue, 0 }  ,fill opacity=1 ] (133,152.09) .. controls (133.05,150.87) and (134.07,149.92) .. (135.29,149.96) .. controls (136.51,150) and (137.46,151.03) .. (137.42,152.25) .. controls (137.37,153.47) and (136.35,154.42) .. (135.13,154.37) .. controls (133.91,154.33) and (132.96,153.3) .. (133,152.09) -- cycle ;
%Shape: Circle [id:dp6450479167540113]
\draw  [fill={rgb, 255:red, 0; green, 0; blue, 0 }  ,fill opacity=1 ] (133,177.26) .. controls (132.95,176.04) and (133.89,175.01) .. (135.11,174.96) .. controls (136.33,174.91) and (137.36,175.85) .. (137.41,177.07) .. controls (137.47,178.29) and (136.52,179.32) .. (135.3,179.37) .. controls (134.09,179.43) and (133.05,178.48) .. (133,177.26) -- cycle ;
%Shape: Circle [id:dp6001863554979303]
\draw  [fill={rgb, 255:red, 0; green, 0; blue, 0 }  ,fill opacity=1 ] (132,230.17) .. controls (132,228.95) and (132.99,227.96) .. (134.21,227.96) .. controls (135.43,227.96) and (136.42,228.95) .. (136.42,230.17) .. controls (136.42,231.39) and (135.43,232.37) .. (134.21,232.37) .. controls (132.99,232.37) and (132,231.39) .. (132,230.17) -- cycle ;
%Rounded Rect [id:dp2690189033599666]
\draw  [dash pattern={on 0.84pt off 2.51pt}] (236,122.87) .. controls (236,116.87) and (240.87,112) .. (246.87,112) -- (279.47,112) .. controls (285.47,112) and (290.33,116.87) .. (290.33,122.87) -- (290.33,235.13) .. controls (290.33,241.13) and (285.47,246) .. (279.47,246) -- (246.87,246) .. controls (240.87,246) and (236,241.13) .. (236,235.13) -- cycle ;
%Shape: Circle [id:dp8491761590857403]
\draw  [fill={rgb, 255:red, 0; green, 0; blue, 0 }  ,fill opacity=1 ] (252,129.17) .. controls (252,127.95) and (252.99,126.96) .. (254.21,126.96) .. controls (255.43,126.96) and (256.42,127.95) .. (256.42,129.17) .. controls (256.42,130.39) and (255.43,131.37) .. (254.21,131.37) .. controls (252.99,131.37) and (252,130.39) .. (252,129.17) -- cycle ;
%Shape: Circle [id:dp2771811314322483]
\draw  [fill={rgb, 255:red, 0; green, 0; blue, 0 }  ,fill opacity=1 ] (252,152.09) .. controls (252.05,150.87) and (253.07,149.92) .. (254.29,149.96) .. controls (255.51,150) and (256.46,151.03) .. (256.42,152.25) .. controls (256.37,153.47) and (255.35,154.42) .. (254.13,154.37) .. controls (252.91,154.33) and (251.96,153.3) .. (252,152.09) -- cycle ;
%Shape: Circle [id:dp8163775575387513]
\draw  [fill={rgb, 255:red, 0; green, 0; blue, 0 }  ,fill opacity=1 ] (252,177.26) .. controls (251.95,176.04) and (252.89,175.01) .. (254.11,174.96) .. controls (255.33,174.91) and (256.36,175.85) .. (256.41,177.07) .. controls (256.47,178.29) and (255.52,179.32) .. (254.3,179.37) .. controls (253.09,179.43) and (252.05,178.48) .. (252,177.26) -- cycle ;
%Shape: Circle [id:dp521357783808911]
\draw  [fill={rgb, 255:red, 0; green, 0; blue, 0 }  ,fill opacity=1 ] (251,230.17) .. controls (251,228.95) and (251.99,227.96) .. (253.21,227.96) .. controls (254.43,227.96) and (255.42,228.95) .. (255.42,230.17) .. controls (255.42,231.39) and (254.43,232.37) .. (253.21,232.37) .. controls (251.99,232.37) and (251,231.39) .. (251,230.17) -- cycle ;
%Straight Lines [id:da4731740736177197]
\draw    (134.21,129.17) -- (254.21,129.17) ;
%Straight Lines [id:da4250920199851522]
\draw    (134.21,129.17) -- (254.13,154.37) ;
%Straight Lines [id:da7920539319815854]
\draw    (134.21,129.17) -- (252,177.26) ;
%Straight Lines [id:da889825698759696]
\draw    (134.21,129.17) -- (251,230.17) ;
%Straight Lines [id:da6811270336848512]
\draw    (135.21,152.17) -- (252,129.17) ;
%Straight Lines [id:da4650197575015669]
\draw    (135.21,152.17) -- (254.13,154.37) ;
%Straight Lines [id:da3016183723305732]
\draw    (135.21,152.17) -- (254.3,179.37) ;
%Straight Lines [id:da8395642075989516]
\draw    (135.21,152.17) -- (253.21,232.37) ;
%Straight Lines [id:da19204661452580374]
\draw    (135.21,177.17) -- (254.21,131.37) ;
%Straight Lines [id:da7915300821460036]
\draw    (135.21,177.17) -- (252,152.09) ;
%Straight Lines [id:da4075760152815433]
\draw    (135.21,177.17) -- (254.21,177.17) ;
%Straight Lines [id:da723665330193124]
\draw    (135.21,177.17) -- (251,230.17) ;
%Straight Lines [id:da4944926848972009]
\draw    (134.21,230.17) -- (256.42,129.17) ;
%Straight Lines [id:da09932614201861001]
\draw    (134.21,230.17) -- (252,152.09) ;
%Straight Lines [id:da4690079491750463]
\draw    (134.21,230.17) -- (252,177.26) ;
%Straight Lines [id:da33607353633225445]
\draw    (134.21,230.17) -- (255.42,230.17) ;

% Text Node
\draw (121,95) node [anchor=north west][inner sep=0.75pt]   [align=left] {$a\sigma$};
% Text Node
\draw (106,121) node [anchor=north west][inner sep=0.75pt]   [align=left] {$a_1$};
% Text Node
\draw (106,142) node [anchor=north west][inner sep=0.75pt]   [align=left] {$a_2$};
% Text Node
\draw (106,169) node [anchor=north west][inner sep=0.75pt]   [align=left] {$a_3$};
% Text Node
\draw (106,220) node [anchor=north west][inner sep=0.75pt]   [align=left] {$a_n$};
% Text Node
\draw (257,95) node [anchor=north west][inner sep=0.75pt]   [align=left] {$b\sigma$};
% Text Node
\draw (268,121) node [anchor=north west][inner sep=0.75pt]   [align=left] {$b_1$};
% Text Node
\draw (268,142) node [anchor=north west][inner sep=0.75pt]   [align=left] {$b_2$};
% Text Node
\draw (268,167) node [anchor=north west][inner sep=0.75pt]   [align=left] {$b_3$};
% Text Node
\draw (268,220) node [anchor=north west][inner sep=0.75pt]   [align=left] {$b_m$};
% Text Node
\draw (132,193.4) node [anchor=north west][inner sep=0.75pt]    {$\vdots $};
% Text Node
\draw (250,193.4) node [anchor=north west][inner sep=0.75pt]    {$\vdots $};
\end{tikzpicture}
\end{center}

Further we have:

\begin{them}\label{GammaS0}
Let $S$ be an inverse semigroup without zero and $|S/\sigma| \geq 2$. The zero-divisor graph  $\Gamma(S^0)$ is the  union of complete bipartite graphs $K_{a\sigma, b\sigma}$, where $(a, b) \notin \sigma$ and $a\sigma$ and $b\sigma$ are two partitions of $K_{a\sigma, b\sigma}$, that is,
\[\Gamma(S^0) =\bigcup_{a, b \in S, (a, b )\notin \sigma}K_{a\sigma, b\sigma} \]
\end{them}

\begin{coro}\label{GroupComplete}
Let $S$ be an inverse semigroup without zero. Then  $\Gamma(S^0)$ is a complete graph if and only if $|S|>1$ and $\sigma = 1_S$, where $1_S$ is the identity relation on $S$, that is, $S$ is not a non\mbox{-}trival group.
\end{coro}

\begin{them}
Let $S$ be an inverse semigroup without zero. Then we have:
\begin{enumerate}
\item[\rm (i)] $\mbox{\bf diam}(\Gamma(S^0) = 1$  if and only if  $|S|>1$ and $\sigma = 1_S$;
\item[\rm (ii)] $\mbox{\bf diam}(\Gamma(S^0) = 2$ if and only if  $|S/\sigma| \geq 2$ and $\sigma \neq 1_S$;
\item[\rm (iii)] $\mbox{\bf diam}(\Gamma(S^0) \in \{1, 2\}$ if and only if $|S/\sigma| \geq 2$,
\end{enumerate}
 where $1_S$ is the identity relation on $S$.
\end{them}
\begin{proof}
(i) It is an immediate result of Corollary~\ref{GroupComplete}.

(ii) If $\mathbf{diam} (\Gamma(S^0)) = 2 $ then there exist $a, b \in S$ such that $a$ and $b$ are not adjacent in $\Gamma(S^0)$ but there exists $c \in S$ with both $a, c$ and $b, c$ being adjacent in $\Gamma(S^0)$. By Lemma~\ref{properInv} we obtain that $(a, c) \notin \sigma$, $(b, c) \notin \sigma$ and $(a, b) \in \sigma$. It follows that $|S/\sigma| \geq 2$.  Conversely, if  $|S/\sigma| \geq 2$  then $S$ contains more than one element and also  by Proposition~\ref{VertexOfProIvs} (ii) we have $V(\Gamma(S^0)) = S$. Since $\sigma \neq 1_S$ it follows from (i) that $\mbox{\bf diam}(\Gamma(S^0)) \neq 1$. For arbitrary distinct elements $a, b \in S$, we have either $(a, b) \in \sigma$ or $(a, b) \notin \sigma$. In the former, it follows from $|S/\sigma| \geq 2$  that there exists $c \in S$ such that $(a, c) \notin \sigma$, and then by Lemma~\ref{properInv} we have $a$ and $b$ are adjacent with $c$, respectively. So ${\bf d}(a, b) = 2$. In the latter, $a$ and $b$ are adjacent by Lemma~\ref{properInv}  and so {\bf d}$(a, b) = 1$. Consequently, $\mbox{\bf diam}(\Gamma(S^0)) = 2$. Hence the result holds.

(iii) It follows from Proposition~\ref{VertexOfProIvs}, part (i) and part (ii) that $\mbox{\bf diam}(\Gamma(S^0)) \in \{1, 2\}$ if and only if $|S/\sigma| \geq 2$.
\end{proof}

Since $\sigma$ is the least group congruence on an inverse semigroup $S$, we obtain that:

\begin{coro}
Let $S$ be an inverse semigroup without zero.
\begin{enumerate}
\item[\rm (i)] $\mbox{\bf diam}(\Gamma(S^0) = 1$  if and only if  $S$ is a non-trivial group if and only if $\Gamma(S^0)$ is a complete graph;
\item[\rm (ii)] $\mbox{\bf diam}(\Gamma(S^0) = 2$ if and only if  $|S/\sigma| \geq 2$ and $S$ is not a group;
\end{enumerate}
\end{coro}

\begin{them}
Let $S$ be an inverse semigroup without zero and $|S/\sigma| \geq 2$.
\begin{enumerate}
\item[\rm (i)] $\mbox{\bf gr}(\Gamma(S^0) = \infty$ if and only if $|S/\sigma| = 2$ and at least one of $\sigma$-classes contains only one element of $S$;
\item[\rm (ii)] $\mbox{\bf gr}(\Gamma(S^0) = 4$ if and only if $|S/\sigma| = 2$ and every $\sigma$-class contains at least two elements of $S$;
\item[\rm (iii)] $\mbox{\bf gr}(\Gamma(S^0) = 3$ if and only if $|S/\sigma| \geq 3$.
\end{enumerate}
\end{them}
\begin{proof}
(i) It follows from Proposition~\ref{girtt} that $\mbox{\bf gr}(\Gamma(S^0) = \infty$ if and only if $\Gamma(S^0)$ is a star graph. Again by Theorem~\ref{GammaS0} we have $\Gamma(S^0)$ is a star graph if and only if $|S/\sigma| = 2$ and at least one of $\sigma$-classes contains only one element of $S$. Hence the result holds.

(ii) It follows from Proposition~\ref{girtt} that $\mbox{\bf gr}(\Gamma(S^0) = 4$  if and only if $\Gamma(S^0)$ is a bipartite but not a star graph. Notice that $\Gamma(S^0)$ is a bipartite if and only if $|S/\sigma| = 2$; and also $\Gamma(S^0)$ is  not a star graph if and only if both of two $\sigma$-classes contain more than one element. Hence the result holds.

(iii) Suppose that $|S/\sigma| \geq 3$. Then there exist $a, b, c \in S$ which are not $\sigma$-related, and so any two of them are adjacent by Lemma~\ref{properInv}. Hence $a\,\tikz[baseline] \draw (0,0.13) -- (0.5,0.13);\, b\, \tikz[baseline] \draw (0,0.13) -- (0.5,0.13);\, c \,\tikz[baseline] \draw (0,0.13) -- (0.5,0.13);\, a$ forms a circle. Certainly it is one of the the shortest circles in $\Gamma(S^0)$ and so $\mbox{gr}(\Gamma(S^0) = 3$. Conversely, if $\mbox{gr}(\Gamma(S^0) = 3$ then there exist $a, b, c \in S$ such that they form a circle, say $a\,\tikz[baseline] \draw (0,0.13) -- (0.5,0.13);\, b\, \tikz[baseline] \draw (0,0.13) -- (0.5,0.13);\, c \,\tikz[baseline] \draw (0,0.13) -- (0.5,0.13);\, a$, which follows from Lemma~\ref{properInv} that $a, b, c$ are not $\sigma$-related, and so $|S/\sigma| \geq 3$.
\end{proof}

\section{Graph inverse semigroups}\label{graph inverse semigroup}
In this section, we focus on zero-divisor graphs of graph inverse semigroups. We begin with recalling the definition of graph inverse semigroups.

Given a directed graph $G = (V(G), E(G), {\bf s}, {\bf r})$, the {\it graph inverse semigroup} $I(G)$ of $G$ is the semigroup with zero generated by sets $V(G)$ and $E(G)$, together with a set $\{e^{-1}: e \in E(G)\}$, satisfying the following relations for all $u, v \in V(G)$ and $e, f \in E(G)$,
\begin{enumerate}
\item[(V)] $uv = \delta_{u,v}u$;
\item[(E1)] ${\bf s}(e)e = e{\bf r}(e) = e$;
\item[(E2)] ${\bf r}(e)e^{-1} = e^{-1}{\bf s}(e) = e^{-1}$;
\item[(CK1)] $e^{-1}f = \delta_{e, f}{\bf r}(e)$,
\end{enumerate}
where $\delta$ is the Kronecker delta.

We emphasize that condition (V) implies that $v^2 = v$ for all $v \in V(G)$, that is, vertices of $G$ are idempotents in $I(G)$. Every non-zero element of $I(G)$ can be written uniquely as $pq^{-1}$ for some $p, q \in $Path$(G)$, by (CK1). It is easy to see that for all $pq^{-1}, rs^{-1} \in I(G)$,

$$ (pq^{-1})(rs^{-1})=
	\begin{cases}
		pts^{-1} & \mbox{if}\ r = qt \mbox{\ for\ some\ (possibly\ empty)\ path}\ t\\
        p(st)^{-1} & \mbox{if}\ q = rt \mbox{\ for\ some\ (possibly\ empty)\ path}\ t\\
		0 & \mbox{otherwise},
	\end{cases}$$
Further the set $E(I(G))$ of idempotents of $I(G)$ is
\[E(I(G))= \{pp^{-1}: p \in \mbox{Path}(G)\} \cup \{0\}.\]
It is also easy to verify that $I(G)$ is indeed an inverse semigroup with $(pq^{-1})^{-1} = qp^{-1}$ for all $p, q \in$ Path$(G)$.

Notice that if  $G =(V(G), E(G), \mathbf{s}, \mathbf{r})$  is a trivial graph, that is,  $G$ consists of one isolated vertex $v$, then $I(G) = \{0, v\}$, and so $I(G)$ does not have zero-divisors.
In the following  we assume that  $G$ is not trivial.

Let $G =(V(G), E(G), \mathbf{s}, \mathbf{r})$ be a  non\mbox{-}trivial directed graph and $\Gamma(I(G))$ be the zero\mbox{-}divisor graph of $I(G)$. To describe the diameter and girth of $\Gamma(I(G))$, we need first describe the set of vertices of $\Gamma(I(G))$.

 We will show that all non-zero element of $I(G)$ are zero-divisors and so $V(\Gamma(I(G))) = I(G)^{\times}$. To do this we start with some lemmas.

It is easy to obtain  Lemma~\ref{Ev} so we omit its proof.
\begin{lem}\label{Ev}
For all $v \in V(G)$, we have $E(I(G))v= \{pp^{-1}: p \in Path(G), {\bf{s}}(p)=v\} \cup \{0\}$.
\end{lem}

\begin{lem}\label{vx}
For all $v\in V(G)$ and $p\in {\rm Path}(G)\setminus\{v\}$, we have $\omega(v, p) = \{0\}$ and $\omega(v, p^{-1})=\{0\}$. Further, $\mbox{Path}(G)\cup\mbox{Path}(G)^{-1} \subseteq Z(I(G))^{\times}$.
\end{lem}
\begin{proof}
Suppose that   $v\in V(G)$, $p\in \mbox{Path}(G)\setminus\{v\}$ and $\alpha\alpha^{-1} \in E(I(G))$. From Lemma~\ref{Ev}, $E(I(G))v= \{qq^{-1}: q \in Path(G),{\bf{s}}(q)=v\} \cup \{0\}$.

If $p \in V(G)$ we have $E(I(G))p= \{tt^{-1}: t \in Path(G), {\bf{s}}(t)=p\} \cup \{0\}$. Since $p\ne v$ it follows that $E(I(G))v\cap E(I(G))p=\{0\}$ and so by Proposition \ref{a}, $\omega(v, p)=\{0\}$. In this case, $p^{-1} = p$, and so $\omega(v, p^{-1}) = \{0\}$, as required.

If $p\in \mbox{Path}(G)\setminus V(G)$, we have

\begin{align}
\alpha\alpha^{-1}p=&
\begin{cases}
p     &\mbox{if}\ p=\alpha\xi\ \mbox{for}\ \mbox{some}\ \xi\in \mbox{Path}(G) \\
p\xi\xi^{-1}   & \mbox{if}\ \alpha=p\xi\ \mbox{for}\ \mbox{some}\ \xi\in \mbox{Path}(G) \\
0          & \mbox{otherwise}.\\
\end{cases}\label{2x}
\end{align}
Since $p\in \mbox{Path}(G)\setminus V(G)$ it follows that $p \notin E(I(G))v $ and  also $ p\xi\xi^{-1} \notin E(I(G))v$. Hence $E(I(G))v\cap E(I(G))p=\{0\}$. By Proposition \ref{a}, $\omega(v, p)=\{0\}$. We also have
\begin{align}
\alpha\alpha^{-1}p^{-1}=&
\begin{cases}
\alpha(p\alpha)^{-1}   & \mbox{if}\ {\bf{r}}(p)={\bf{s}}(\alpha) \\
0          &\mbox{otherwise}. \\
\end{cases}\label{2x^{-1}}		
\end{align}
Since $p \notin V(G)$ and ${\bf{r}}(p)={\bf{s}}(\alpha)$ it follows that $\alpha \neq p\alpha$ and so $\alpha(p \alpha)^{-1} \notin E(I(G))v$ which implies that  $E(I(G))v\cap E(I(G))p^{-1}=\{0\}$. By Proposition \ref{a}, $\omega(v, p^{-1})= \{0\}$, as required.
\end{proof}

\begin{prop}\label{complete}
The induced subgraph  with  set $\mbox{Path}(G)\cup\mbox{Path}(G)^{-1}$ of vertices is a complete subgraph of $\Gamma(I(G))$.
\end{prop}
\begin{proof}
It is sufficient to show that any two distinct  elements in  $\mbox{Path}(G)\cup\mbox{Path}(G)^{-1}$ are adjacent in $\Gamma(I(G))$. It follows from  Lemma \ref{vx} and Proposition \ref{a} that for all $v\in V(G)$, $v$ is adjacent with any other element in $\mbox{Path}(G)\cup\mbox{Path}(G)^{-1}$. Next we only consider adjacency relations on $(\mbox{Path}(G)\cup\mbox{Path}(G)^{-1})\setminus V(G)$.  Let $p, q\in \mbox{Path}(G)\setminus V(G)$ with $p\ne q$.
	
Case a.  We show that $p$ and $q$ are adjacent in $\Gamma(I(G))$. According to ~(\ref{2x}) we get
\begin{equation}
    E(I(G))p= \{p\} \cup \{p\xi\xi^{-1}: {\bf{r}}(p) = {\bf{s}}(\xi), \xi \in \mbox{Path}(G)\} \cup \{0\}
\end{equation}
and
\begin{equation}
E(I(G))q= \{q\} \cup \{q\eta\eta^{-1}: {\bf{r}}(q) = {\bf{s}}(\eta), \eta \in \mbox{Path}(G)\} \cup \{0\}.
\end{equation}
Notice that $p = q\eta\eta^{-1}$ only when $\eta \in V(G)$ and $p = q$. It is a contradiction with $p \neq q$. So $p \neq q\eta\eta^{-1}$. Similarly, $q \neq p\xi\xi^{-1}$. If $p\xi\xi^{-1} = q\eta\eta^{-1}$  we must have $\xi = \eta$ and $p = q$ which is a contradiction with  $p \neq q$. Hence $E(I(G))p \cap E(I(G))q = \{0\}$ and so by Proposition \ref{a}, $p$ and $q$ are adjacent in $\Gamma(I(G))$.
	
Case b. We now prove that $p$ and $q^{-1}$ are adjacent  in $\Gamma(I(G))$. According to  ~(\ref{2x^{-1}}) we get
\begin{equation}
E(I(G))q^{-1}= \{\mu(q\mu)^{-1}: {\bf{r}}(q) = {\bf{s}}(\mu), \mu \in \mbox{Path}(G)\} \cup \{0\}.
\end{equation}
Compare (3) and (5). Notice that $p \neq \mu(q\mu)^{-1}$ as $p,q \in \mbox{Path}(G) \backslash V(G)$.  If $p\xi\xi^{-1} = \mu(q\mu)^{-1}$, we get that $p\xi = \mu$ and $\xi = q\mu$, which implies that $pq\mu = \mu$ which is contradictory to $p, q \notin V(G)$. Thus  $p\xi\xi^{-1} \ne \mu(q\mu)^{-1}$. Hence  $E(I(G))p \cap E(I(G))q^{-1} = \{0\}$ and so by Proposition \ref{a}, $p$ and $q^{-1}$ are adjacent in $\Gamma(I(G))$.

Case c.  We now show that $p^{-1}$ and $q^{-1}$ are adjacent in $\Gamma(I(G))$. According to ~(\ref{2x^{-1}}), we get
\begin{equation}
 E(I(G))p^{-1}= \{\nu(p\nu)^{-1}: {\bf{r}}(p) = {\bf{s}}(\nu), \nu \in \mbox{Path}(G)\} \cup \{0\}
\end{equation}
If $\nu(p\nu)^{-1} = \mu(q\mu)^{-1}$, we get that $\nu = \mu$ and $p\nu = q\mu$, and so $p = q$ which is contradictory to $p \ne q$.  Hence  $E(I(G))p^{-1} \cap E(I(G))q^{-1} = \{0\}$ and so by Proposition \ref{a}, $p^{-1}$ and $q^{-1}$ are adjacent in $\Gamma(I(G))$.
\end{proof}

\begin{lem}\label{XXXX}
 For all $p,q,\xi\in \mbox{Path}(G)$ with ${\bf{r}}(p)= {\bf{r}}(q)= {\bf{s}}(\xi)$ and $pq^{-1}\ne p\xi$, we have $\omega(pq^{-1}, p\xi) = \{0\}$. Further $pq^{-1} \in Z(I(G))^{\times}$.
\end{lem}
\begin{proof}
Suppose that  $p,q,\xi\in \mbox{Path}(G)$ with ${\bf{r}}(p)= {\bf{r}}(q)= {\bf{s}}(\xi)$ and $pq^{-1}\ne p\xi$. If $p \in V(G)$ or $q \in V(G)$, we get $pq^{-1} \in \mbox{Path}(G)\cup\mbox{Path}(G)^{-1}$, and so by Proposition \ref{complete} $pq^{-1}$ and $p\xi$ are adjacent in $\Gamma(I(G))$.
Now we assume that $p, q \notin V(G)$. For all $\alpha\alpha^{-1},\beta\beta^{-1}\in E(I(G))$, we have
	$$
	\alpha\alpha^{-1}pq^{-1} =
	\begin{cases}
		pq^{-1}       &\mbox{if}\ p=\alpha\eta\ \mbox{for}\ \mbox{some}\ \eta\in \mbox{Path}(G) \\
		(p\eta)(q\eta)^{-1}   & \mbox{if}\ \alpha=p\eta\ \mbox{for}\ \mbox{some}\ \eta\in \mbox{Path}(G) \\
		0          & \mbox{otherwise} \\
	\end{cases}\\
	$$
and
	$$
	\beta\beta^{-1}p\xi =
	\begin{cases}
		p\xi       &\mbox{if}\ p\xi=\beta\gamma\ \mbox{for}\ \mbox{some}\ \gamma\in \mbox{Path}(G) \\
		p\xi \gamma\gamma^{-1}   & \mbox{if}\ \beta=p\xi\gamma\ \mbox{for}\ \mbox{some}\ \gamma\in \mbox{Path}(G) \\
		0          & \mbox{otherwise}. \\
	\end{cases}\\
	$$
So
$$E(I(G))pq^{-1} = \{(p\eta)(q\eta)^{-1}: \eta \in \mbox{Path}(G)\} \cup \{0\}$$
 and 	
$$E(I(G))p\xi = \{p\xi\gamma\gamma^{-1}: \gamma \in \mbox{Path}(G)\} \cup \{0\}.$$
If $(p\eta)(q\eta)^{-1} = p\xi\gamma\gamma^{-1}$ we obtain that $p\eta = p\xi\gamma$ and $q\eta = \gamma$, and so $p\eta = p\xi q\eta$, which implies that $p = p\xi q$. As $p \notin V(G)$ we must have $\xi, q \in V(G)$, which is contradictory to $q \notin V(G)$. Hence $E(I(G))pq^{-1} \cap E(I(G))p\xi = \{0\}$. By Proposition \ref{a} $\omega(pq^{-1}, p\xi) = \{0\}$, and so $pq^{-1} \in Z(I(G))^{\times}$, as required.
\end{proof}

Certainly we have $Z(I(G))^{\times} \subseteq I(G)^{\times}$. Conversely $I(G)^{\times} \subseteq Z(I(G))^{\times}$ is obtained by Lemma~\ref{XXXX}. Then we have :

\begin{prop}\label{vertex set}
Every non-zero element of $I(G)$ is a vertex of $\Gamma(I(G))$, that is, $V(\Gamma(I(G))) = I(G)^{\times}$.
\end{prop}

\begin{lem}\label{lem1}
 For all distinct vertices $pq^{-1}, rs^{-1}\in V(\Gamma(I(G)))$, $pq^{-1}$ and $rs^{-1}$ are not adjacent in $\Gamma(I(G))$ if and only if there exist $\xi,\eta\in \mbox{Path}(G)$ such that $p\xi=r\eta$ and $q\xi=s\eta$, where ${\bf{r}}(p)={\bf{r}}(q)={\bf{s}}(\xi)$ and ${\bf{r}}(r)={\bf{r}}(s)={\bf{s}}(\eta)$.
	\end{lem}
\begin{proof}
Necessity. Suppose that $pq^{-1}$ and $rs^{-1}$ are not adjacent in $\Gamma(I(G))$. By Proposition \ref{a} we have $E(I(G))pq^{-1}\cap E(I(G))rs^{-1}-\{0\}\ne \emptyset$ and so there exist $\alpha\alpha^{-1},\beta\beta^{-1}\in E(I(G))$ such that $\alpha\alpha^{-1}pq^{-1}=\beta\beta^{-1}rs^{-1}\ne 0$, where
$$
\alpha\alpha^{-1}pq^{-1} =
\begin{cases}
	pq^{-1}       &\mbox{if}\ p=\alpha x\ \mbox{for}\ \mbox{some}\ x\in \mbox{Path}(G) \\
	(px)(qx)^{-1}   & \mbox{if}\ \alpha=px\ \mbox{for}\ \mbox{some}\ x\in \mbox{Path}(G) \\
	0          & \mbox{otherwise}, \\
\end{cases}\\
$$
$$\beta\beta^{-1}rs^{-1} =
\begin{cases}
	rs^{-1}       &\mbox{if}\ r=\beta y\ \mbox{for}\ \mbox{some}\ y\in \mbox{Path}(G) \\
	(ry)(sy)^{-1}   & \mbox{if}\ \beta=ry\ \mbox{for}\ \mbox{some}\ y\in \mbox{Path}(G) \\
	0          & \mbox{otherwise}. \\
\end{cases}$$
If $\alpha\alpha^{-1}pq^{-1}=\beta\beta^{-1}rs^{-1}\ne 0$  there exist four cases to discuss.

Case a. If $p=\alpha x$ and $r=\beta y$ for some $x, y\in \mbox{Path}(G)$, we have
$$pq^{-1}=\alpha\alpha^{-1}pq^{-1}=\beta\beta^{-1}rs^{-1} = rs^{-1}.$$
It follows that $p=r$ and $q=s$. Take $\xi = \eta = {\bf{r}}(p)$ as ${\bf{r}}(p) = {\bf{r}}(q)$, and so the result is as required.

Case b. If $p=\alpha x$  and $\beta=ry$ for some $x, y\in \mbox{Path}(G)$, we have
$$pq^{-1}=\alpha\alpha^{-1}pq^{-1}=\beta\beta^{-1}rs^{-1} = (ry)(sy)^{-1}.$$
It follows that $p=ry$ and $q=sy$. Take $\xi = {\bf{r}}(p)$ and $\eta = y$ as ${\bf{r}}(p) = {\bf{r}}(q)$, and so the result is as required.

Case c. If $\alpha = px$ and  $r=\beta y$ for some $x, y\in \mbox{Path}(G)$, we have
$$(px)(qx)^{-1}=\alpha\alpha^{-1}pq^{-1}=\beta\beta^{-1}rs^{-1} = rs^{-1}.$$
It follows that $px = r$ and $qx = s$. Take $\xi = x$ and $\eta = {\bf{r}}(r)$ as ${\bf{r}}(r) = {\bf{r}}(s)$, and so the result is as required.

Case d. If $\alpha = px$ and $\beta=ry$ for some $x, y\in \mbox{Path}(G)$, we have
$$(px)(qx)^{-1}=\alpha\alpha^{-1}pq^{-1}=\beta\beta^{-1}rs^{-1} = (ry)(sy)^{-1}.$$
It follows that $px = ry$ and $qx = sy$. Take $\xi = x$ and $\eta = y$ and then the result is  as required.

Sufficiency. Suppose that there exist $\xi,\eta\in \mbox{Path}(G)$ such that $p\xi=r\eta$ and $q\xi=s\eta$, where ${\bf{r}}(p)={\bf{r}}(q)={\bf{s}}(\xi)$ and ${\bf{r}}(r)={\bf{r}}(s)={\bf{s}}(\eta)$. Set $\alpha=p\xi=r\eta$.
We have $\alpha\alpha^{-1}pq^{-1}=(p\xi)(q\xi)^{-1}=(r\eta)(s\eta)^{-1}= \alpha\alpha^{-1}rs^{-1} \neq 0$. By Proposition \ref{a} $pq^{-1}$ is not adjacent to $rs^{-1}$ in $\Gamma(I(G))$.
\end{proof}

In the following we characterize the diameter and girth of $\Gamma(I(G))$.

\begin{them}\label{them1}
Let $G = (V(G), E(G), \mathbf{s}, \mathbf{r})$ be a
non\mbox{-}trivial directed graph. Then $\mbox{\bf diam}(\Gamma(I(G)))$ $\in \{1, 2\}$.
\end{them}
\begin{proof}By Proposition \ref{vertex set}, $V(\Gamma(I(G)))=I(G)^{\times}$. Suppose that  $pq^{-1}$ and $rs^{-1}$ are two distinct vertices of $\Gamma(I(G))$. If $pq^{-1}$ is adjacent to $rs^{-1}$ in $\Gamma(I(G)$, then we have ${\bf d}(pq^{-1},rs^{-1})=1$. Now we suppose that $pq^{-1}$ and $rs^{-1}$ are not adjacent in $\Gamma(I(G)$.  Then it follows from Lemma \ref{lem1} that there exist $\xi, \eta\in \mbox{Path}(G)$ such that $p\xi=r\eta$ and $q\xi=s\eta$, where ${\bf{r}}(p)={\bf{r}}(q)={\bf{s}}(\xi)$ and ${\bf{r}}(r)={\bf{r}}(s)={\bf{s}}(\eta)$. According to Lemma  \ref{XXXX} if $pq^{-1}\ne p\xi$ and $rs^{-1}\ne r\eta$, then $\omega(pq^{-1}, p\xi)=0$ and $\omega(rs^{-1}, r\eta)=0$, that is, $pq^{-1}$ is adjacent to $p\xi$ and $rs^{-1}$ is adjacent to $r\eta$ in $\Gamma(I(G))$. Together with $p\xi = r\eta$,  there exists a path $pq^{-1}\,\tikz[baseline] \draw (0,0.13) -- (0.5,0.13);\, p\xi \,\tikz[baseline] \draw (0,0.13) -- (0.5,0.13);\, rs^{-1}$  in $\Gamma(I(G))$. It follows that ${\bf d}(pq^{-1},rs^{-1})=2$. So it is sufficient to  show that $pq^{-1}\ne p\xi$ and $rs^{-1}\ne r\eta$. Since $pq^{-1} \neq rs^{-1}$ and $p\xi=r\eta$ it follows that $pq^{-1}=p\xi$ and $rs^{-1}= r\eta$ can not occur simultaneously.
If $pq^{-1}= p\xi$ and $rs^{-1}\ne r\eta$, then by Lemma \ref{XXXX}, $\omega(r\eta, rs^{-1}) = 0$, that is, $r\eta$ is adjacent to $rs^{-1}$ in $\Gamma(I(G)$. Together with $r\eta=p\xi=pq^{-1}$, we get $pq^{-1}$ is adjacent to $rs^{-1}$ in $\Gamma(I(G)$, which is contradictory to the assumption $pq^{-1}$ and $rs^{-1}$ being not adjacent. Hence $pq^{-1}= p\xi$ and $rs^{-1}\ne r\eta$ can not occur simultaneously. Similarly,  $pq^{-1}\ne p\xi$ and $rs^{-1}=r\eta$ can not occur simultaneously. Hence, $pq^{-1}\ne p\xi$ and $rs^{-1}\ne r\eta$, as required. Consequently,  $\Gamma(I(G))$ is a graph with $\mbox{\bf diam}(\Gamma(S))\in\{1,2\}$.
\end{proof}

\begin{them}\label{them3}
Let $G = (V(G), E(G), \mathbf{s}, \mathbf{r})$ be a
non\mbox{-}trivial directed graph. Then  we have
\begin{enumerate}
\item[\rm (i)]  $\mbox{\bf diam}(\Gamma(I(G)))=1$ if and only if $I(G)^{\times}=\mbox{\rm Min}(I(G)^{\times})$, if and only if $G$ is a null graph with $|V(G)|\geq 2$;
\item[\rm (ii)]  $\mbox{\bf diam}(\Gamma(I(G))) = 2$ if and only if $G$ is not a null graph, if and only if there exist distinct elements $pq^{-1}, rs^{-1}\in I(G)^{\times}$ such that $pq^{-1}$ and $rs^{-1}$ are comparable.
\end{enumerate}
\end{them}
\begin{proof}
(i) By Proposition \ref{vertex set} we get that  $V(\Gamma(I(G))) = I(G)^{\times}$. It follows from Proposition~\ref{diamP} that $\mbox{\bf diam}(\Gamma(I(G)))=1$ if and only if $I(G)^{\times}=\mbox{Min}(I(G)^{\times})$.

Now we show that $\mbox{\bf diam}(\Gamma(I(G)))=1$  if and only if $G$ is a null graph with $|V(G)|\geq 2$. Let $G$ be a null  graph with $|V(G)|\geq 2$, that is, graph $G$ consists of $n$ isolated vertices, where $n \ge 2$. Then $I(G)=\{0,v_1,\ldots,v_n\}$. By Proposition \ref{vertex set} we get $V(\Gamma(I(G)))=\{v_1,\ldots,v_n\}$.  Thus $\Gamma(I(G))$ is complete with $n$ vertices by Proposition\ref{complete} and so $\mbox{\bf diam}(\Gamma(S))=1$.
	
Conversely, suppose that $\mbox{\bf diam}(\Gamma(I(G)))=1$ and let $|V(G)|=n$. There exist two case: either $n = 1$ or $n \geq 2$. If $n = 1$, we denote the unique vertex  by $v$. Then there must exist an edge $e$ with $\mathbf{s}(e) = \mathbf{r}(e) = v$, otherwise $\Gamma(I(G))$ does not exist. Further we have $I(G)=\{0, v, e, e^{-1}, ee^{-1}\}$. Since $vee^{-1} = ee^{-1}$ we get that $ee^{-1} \leq v$ and so $\omega(v, ee^{-1}) = \{0,ee^{-1}\}$, that is, $v$ and $ee^{-1}$ are not adjacent in $\Gamma(I(G))$. Then ${\bf{d}}(v, ee^{-1}) > 1$. It follows that $\mbox{\bf diam}(\Gamma(I(G)))>1$, a contradiction. Hence $n \geq 2$.

Suppose that  $n \geq 2$ and $|E(G)|\neq 0$. Then there exists an edge $e\in E(G)$ with $v_i={\bf{s}}(e)$ and $v_j={\bf{r}}(e)$. So $0 \neq ee^{-1} \in E(I(G))$. We also have $ee^{-1} \in E(I(G))ee^{-1}$ and $ee^{-1} =ee^{-1}v_i \in E(I(G))v_i$. Thus $0\ne ee^{-1}\in E(I(G))v_i\cap E(I(G))ee^{-1}$.  According to Proposition \ref{a}, $v_i$ is not adjacent to $ee^{-1}$ in $\Gamma(I(G))$. So {\bf d}$(v_i, ee^{-1}) = 2$ by Theorem~\ref{them1}, which is a contradiction to $\mbox{\bf diam}(\Gamma(I(G)))=1$. Hence $|E(G)|=0$, that is, $G$ is a null graph with $|V(G)|\geq 2$.

(ii) It is easy to see that $\mbox{\bf diam}(\Gamma(I(G))) = 2$ if and only if $G$ is not a null graph by Theorem~\ref{them1} and part (i).

Now suppose that $\mbox{\bf diam}(\Gamma(I(G)))=2$. Then there exist $pq^{-1}$, $\mu\nu^{-1}\in I(G)^{\times}$ with {\bf d}$(pq^{-1}, \mu\nu^{-1})=2$ which implies $\omega(pq^{-1}, \mu\nu^{-1}) \ne\{0\}$ by Proposition \ref{a}. So there exists $ rs^{-1}\in I(G)^{\times}$ such that $rs^{-1}\in \omega(pq^{-1}, \mu\nu^{-1})$ which indicates $rs^{-1}\leq pq^{-1}$, as required.
	
Conversely, suppose that there exist distinct elements  $pq^{-1},rs^{-1}\in I(G)^{\times}$ are such that  $pq^{-1}\leq rs^{-1}$ . Then we have $\omega(pq^{-1}, rs^{-1})\ne\{0\}$, which implies that {\bf d}$(pq^{-1},rs^{-1})\ge 2$ by Proposition \ref{a}. It leads that $\mbox{\bf diam}(\Gamma(I(G))) \geq2$.  Again by Theorem \ref{them1} we get  $\mbox{\bf diam}(\Gamma(I(G)))=2$.
\end{proof}

\begin{them}\label{them2}
Let $G = (V(G), E(G), \mathbf{s}, \mathbf{r})$ be a
non\mbox{-}trivial directed graph. The following statements hold:
\begin{enumerate}
\item[\rm (i)]  $\mbox{\bf gr}(\Gamma(I(G)))\in \{3,\infty\}$;
\item[\rm (ii)]  $\mbox{\bf gr}(\Gamma(I(G))) = \infty$ if and only if $|V(G)|=2$ and  $|E(G)|=0$;
\item[\rm (iii)]  $\mbox{\bf gr}(\Gamma(I(G))) = 3$ if and only if either $|E(G)|\neq 0$, or $|V(G)|\ge 3$ and $|E(G)|=0$.
\end{enumerate}
\end{them}
\begin{proof}
Part (iii) can be obtained  by part (i) and part (ii) so it is sufficient to show part (i) and (ii).

(i) Let $G = (V(G), E(G), \mathbf{s}, \mathbf{r})$ be a non-trivial  directed graph. We discuss it by two cases: either $|E(G)| \neq 0$ or $|E(G)| = 0$.

If $|E(G)|\ne 0$, then there exists $e \in E(G)$. Assume that  ${\bf{r}}(e)=v\in V(G)$. Put $U=\{e, e^{-1}, v\}$. Clearly, $U \subseteq \mbox{Path}(G) \cup \mbox{path}(G)^{-1}$. It follows from Proposition~\ref{complete} that  there exists a cycle of length of $3$ in $\Gamma(I(G))$: $e\,\tikz[baseline] \draw (0,0.13) -- (0.5,0.13);\, e^{-1} \,\tikz[baseline] \draw (0,0.13) -- (0.5,0.13);\, v \,\tikz[baseline] \draw (0,0.13) -- (0.5,0.13);\, e$. It follows from Proposition~\ref{girtt} that $\mbox{gr}(\Gamma(I(G)))=3$.

If $|E(G)|=0$ and $V(G)=n$  with $n\ge 2$. Then we have $I(G)=\{0,v_1,\ldots,v_n\}$ and $V(\Gamma(I(G)))=\{v_1,\ldots,v_n\}$. According to Proposition~\ref{complete}, $\Gamma(I(G))$ is a complete graph with $n$ vertices. So by Proposition~\ref{girtt}(ii). we have $\mbox{gr}(\Gamma(I(G))) = \infty$ when $n=2$, and by Proposition 3.7(iv) $\mbox{\bf gr}(\Gamma(I(G))) =3$ when $n \ge 3$. Hence $\mbox{\bf gr}(\Gamma(I(G)))\in\{3,\infty\}$.
	
(ii) It follows from the proof of part (i) that the sufficiency holds. Now we show that if $\mbox{gr}(\Gamma(I(G)))=\infty$ then $|V(G)|=2$ and  $|E(G)|=0$. Suppose that $|E(G)|\ne 0$. Then by the proof of part (i) we  get $\mbox{\bf gr}(\Gamma(I(G)))=3$, which is a contradiction to the assumption that $\mbox{gr}(\Gamma(I(G)))=\infty$. Hence $|E(G)|=0$.

 Suppose that $|V(G)| > 2$. Then there exist at least three distinct vertices $v_1, v_2, v_3\in V(G)$. By Proposition~\ref{complete}, we have a cycle of length 3: $v_1\,\tikz[baseline] \draw (0,0.13) -- (0.5,0.13);\, v_2 \,\tikz[baseline] \draw (0,0.13) -- (0.5,0.13);\, v_3 \,\tikz[baseline] \draw (0,0.13) -- (0.5,0.13);\, v_1$. It follows from Proposition~\ref{girtt}(iv) that  $\mbox{\bf gr}(\Gamma(I(G)))=3$ which is a contradiction to the assumption that $\mbox{gr}(\Gamma(I(G)))=\infty$. Hence $|V(G)|=2$.
\end{proof}

In view of Theorem \ref{them1}, Theorem \ref{them3} and Theorem \ref{them2}, we have:
\begin{coro}\label{Cor5.10}
Let $G = (V(G), E(G), \mathbf{s}, \mathbf{r})$ be a non\mbox{-}trivial directed graph, there exist three situations about the diameter and girth of the zero-divisor graph $\Gamma$ of $I(G)$:
\begin{enumerate}
\item[\rm(i)]$(\mbox{\bf diam}(\Gamma),(\mbox{\bf gr}(\Gamma))=(1,\infty)$ if and only if $|E(G)|=0$ and $|V(G)|=2$;
\item[\rm(ii)]$(\mbox{\bf diam}(\Gamma),(\mbox{\bf gr}(\Gamma))=(1,3)$ if and only if  $|E(G)|=0$ and $|V(G)|\ge3$;
\item[\rm(iii)]$(\mbox{\bf diam}(\Gamma),(\mbox{\bf gr}(\Gamma))=(2,3)$ if and only if  $|E(G)|\ne0$.
\end{enumerate}
\end{coro}

At the end of this section we  give examples of these three cases in Corollary~\ref{Cor5.10}.

\begin{ex}\label{3ex}
Graphs (A), (B) and (C) are the corresponding zero-divisor graphs of $I(G_1)$, $I(G_2)$ and $I(G_3)$ where $G_1$, $G_2$ and $G_3$ are the following graphs (a), (b) and (c), respectively.
\end{ex}
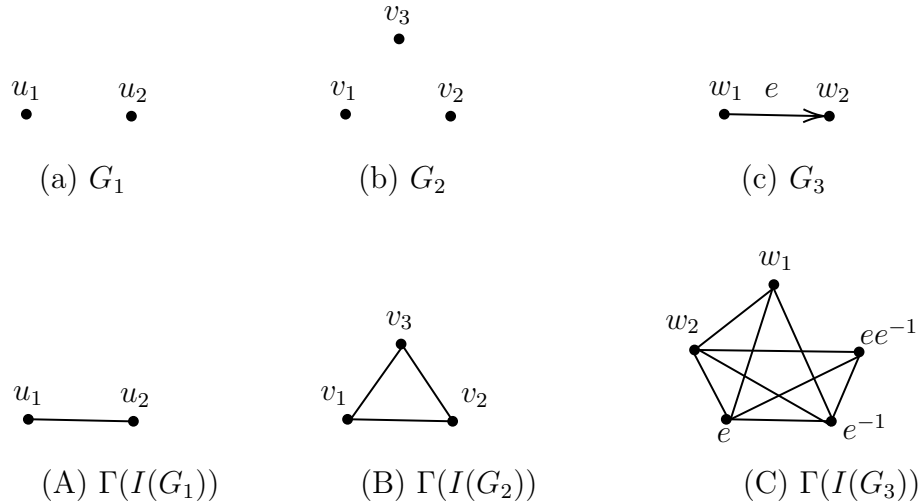
\tikzset{every picture/.style={line width=0.75pt}} %set default line width to 0.75pt
\begin{center}
\begin{tikzpicture}[x=0.75pt,y=0.75pt,yscale=-1,xscale=1]

%Shape: Circle [id:dp963618413238648]
\draw  [fill={rgb, 255:red, 0; green, 0; blue, 0 }  ,fill opacity=1 ] (150,92.17) .. controls (150,90.95) and (150.99,89.96) .. (152.21,89.96) .. controls (153.43,89.96) and (154.42,90.95) .. (154.42,92.17) .. controls (154.42,93.39) and (153.43,94.37) .. (152.21,94.37) .. controls (150.99,94.37) and (150,93.39) .. (150,92.17) -- cycle ;
%Shape: Circle [id:dp03087323152916066]
\draw  [fill={rgb, 255:red, 0; green, 0; blue, 0 }  ,fill opacity=1 ] (203,93.17) .. controls (203,91.95) and (203.99,90.96) .. (205.21,90.96) .. controls (206.43,90.96) and (207.42,91.95) .. (207.42,93.17) .. controls (207.42,94.39) and (206.43,95.37) .. (205.21,95.37) .. controls (203.99,95.37) and (203,94.39) .. (203,93.17) -- cycle ;
%Shape: Circle [id:dp05648700942816265]
\draw  [fill={rgb, 255:red, 0; green, 0; blue, 0 }  ,fill opacity=1 ] (502,92.17) .. controls (502,90.95) and (502.99,89.96) .. (504.21,89.96) .. controls (505.43,89.96) and (506.42,90.95) .. (506.42,92.17) .. controls (506.42,93.39) and (505.43,94.37) .. (504.21,94.37) .. controls (502.99,94.37) and (502,93.39) .. (502,92.17) -- cycle ;
%Shape: Circle [id:dp9197137681096088]
\draw  [fill={rgb, 255:red, 0; green, 0; blue, 0 }  ,fill opacity=1 ] (555,93.17) .. controls (555,91.95) and (555.99,90.96) .. (557.21,90.96) .. controls (558.43,90.96) and (559.42,91.95) .. (559.42,93.17) .. controls (559.42,94.39) and (558.43,95.37) .. (557.21,95.37) .. controls (555.99,95.37) and (555,94.39) .. (555,93.17) -- cycle ;
%Straight Lines [id:da7008279992493405]
\draw    (506.42,92.17) -- (553,93.13) ;
\draw [shift={(555,93.17)}, rotate = 181.18] [color={rgb, 255:red, 0; green, 0; blue, 0 }  ][line width=0.75]    (10.93,-3.29) .. controls (6.95,-1.4) and (3.31,-0.3) .. (0,0) .. controls (3.31,0.3) and (6.95,1.4) .. (10.93,3.29)   ;
%Shape: Circle [id:dp13302159130569624]
\draw  [fill={rgb, 255:red, 0; green, 0; blue, 0 }  ,fill opacity=1 ] (311,92.17) .. controls (311,90.95) and (311.99,89.96) .. (313.21,89.96) .. controls (314.43,89.96) and (315.42,90.95) .. (315.42,92.17) .. controls (315.42,93.39) and (314.43,94.37) .. (313.21,94.37) .. controls (311.99,94.37) and (311,93.39) .. (311,92.17) -- cycle ;
%Shape: Circle [id:dp36406111458845203]
\draw  [fill={rgb, 255:red, 0; green, 0; blue, 0 }  ,fill opacity=1 ] (364,93.17) .. controls (364,91.95) and (364.99,90.96) .. (366.21,90.96) .. controls (367.43,90.96) and (368.42,91.95) .. (368.42,93.17) .. controls (368.42,94.39) and (367.43,95.37) .. (366.21,95.37) .. controls (364.99,95.37) and (364,94.39) .. (364,93.17) -- cycle ;
%Shape: Circle [id:dp5405305179643978]
\draw  [fill={rgb, 255:red, 0; green, 0; blue, 0 }  ,fill opacity=1 ] (338,54.17) .. controls (338,52.95) and (338.99,51.96) .. (340.21,51.96) .. controls (341.43,51.96) and (342.42,52.95) .. (342.42,54.17) .. controls (342.42,55.39) and (341.43,56.37) .. (340.21,56.37) .. controls (338.99,56.37) and (338,55.39) .. (338,54.17) -- cycle ;
%Shape: Circle [id:dp506356720617001]
\draw  [fill={rgb, 255:red, 0; green, 0; blue, 0 }  ,fill opacity=1 ] (151,246.17) .. controls (151,244.95) and (151.99,243.96) .. (153.21,243.96) .. controls (154.43,243.96) and (155.42,244.95) .. (155.42,246.17) .. controls (155.42,247.39) and (154.43,248.37) .. (153.21,248.37) .. controls (151.99,248.37) and (151,247.39) .. (151,246.17) -- cycle ;
%Shape: Circle [id:dp9751979220440055]
\draw  [fill={rgb, 255:red, 0; green, 0; blue, 0 }  ,fill opacity=1 ] (204,247.17) .. controls (204,245.95) and (204.99,244.96) .. (206.21,244.96) .. controls (207.43,244.96) and (208.42,245.95) .. (208.42,247.17) .. controls (208.42,248.39) and (207.43,249.37) .. (206.21,249.37) .. controls (204.99,249.37) and (204,248.39) .. (204,247.17) -- cycle ;
%Shape: Circle [id:dp8598767091295321]
\draw  [fill={rgb, 255:red, 0; green, 0; blue, 0 }  ,fill opacity=1 ] (503,246.17) .. controls (503,244.95) and (503.99,243.96) .. (505.21,243.96) .. controls (506.43,243.96) and (507.42,244.95) .. (507.42,246.17) .. controls (507.42,247.39) and (506.43,248.37) .. (505.21,248.37) .. controls (503.99,248.37) and (503,247.39) .. (503,246.17) -- cycle ;
%Shape: Circle [id:dp23605717434906826]
\draw  [fill={rgb, 255:red, 0; green, 0; blue, 0 }  ,fill opacity=1 ] (556,247.17) .. controls (556,245.95) and (556.99,244.96) .. (558.21,244.96) .. controls (559.43,244.96) and (560.42,245.95) .. (560.42,247.17) .. controls (560.42,248.39) and (559.43,249.37) .. (558.21,249.37) .. controls (556.99,249.37) and (556,248.39) .. (556,247.17) -- cycle ;
%Shape: Circle [id:dp3656040157998113]
\draw  [fill={rgb, 255:red, 0; green, 0; blue, 0 }  ,fill opacity=1 ] (312,246.17) .. controls (312,244.95) and (312.99,243.96) .. (314.21,243.96) .. controls (315.43,243.96) and (316.42,244.95) .. (316.42,246.17) .. controls (316.42,247.39) and (315.43,248.37) .. (314.21,248.37) .. controls (312.99,248.37) and (312,247.39) .. (312,246.17) -- cycle ;
%Shape: Circle [id:dp7429261541926264]
\draw  [fill={rgb, 255:red, 0; green, 0; blue, 0 }  ,fill opacity=1 ] (365,247.17) .. controls (365,245.95) and (365.99,244.96) .. (367.21,244.96) .. controls (368.43,244.96) and (369.42,245.95) .. (369.42,247.17) .. controls (369.42,248.39) and (368.43,249.37) .. (367.21,249.37) .. controls (365.99,249.37) and (365,248.39) .. (365,247.17) -- cycle ;
%Shape: Circle [id:dp7978852018315239]
\draw  [fill={rgb, 255:red, 0; green, 0; blue, 0 }  ,fill opacity=1 ] (339,208.17) .. controls (339,206.95) and (339.99,205.96) .. (341.21,205.96) .. controls (342.43,205.96) and (343.42,206.95) .. (343.42,208.17) .. controls (343.42,209.39) and (342.43,210.37) .. (341.21,210.37) .. controls (339.99,210.37) and (339,209.39) .. (339,208.17) -- cycle ;
%Straight Lines [id:da8162684563303702]
\draw    (155.42,246.17) -- (204,247.17) ;
%Straight Lines [id:da11445023988520941]
\draw    (316.42,246.17) -- (367.21,247.17) ;
%Straight Lines [id:da5596677364117437]
\draw    (314.21,248.37) -- (341.21,210.37) ;
%Straight Lines [id:da26510920127954174]
\draw    (341.42,208.17) -- (367.42,247.17) ;
%\draw    (343.42,208.17) -- (369.42,247.17) ;
%Shape: Circle [id:dp4820319213991353]
\draw  [fill={rgb, 255:red, 0; green, 0; blue, 0 }  ,fill opacity=1 ] (487,211.17) .. controls (487,209.95) and (487.99,208.96) .. (489.21,208.96) .. controls (490.43,208.96) and (491.42,209.95) .. (491.42,211.17) .. controls (491.42,212.39) and (490.43,213.37) .. (489.21,213.37) .. controls (487.99,213.37) and (487,212.39) .. (487,211.17) -- cycle ;
%Shape: Circle [id:dp8931339295348826]
\draw  [fill={rgb, 255:red, 0; green, 0; blue, 0 }  ,fill opacity=1 ] (570,212.17) .. controls (570,210.95) and (570.99,209.96) .. (572.21,209.96) .. controls (573.43,209.96) and (574.42,210.95) .. (574.42,212.17) .. controls (574.42,213.39) and (573.43,214.37) .. (572.21,214.37) .. controls (570.99,214.37) and (570,213.39) .. (570,212.17) -- cycle ;
%Shape: Circle [id:dp5513342426882331]
\draw  [fill={rgb, 255:red, 0; green, 0; blue, 0 }  ,fill opacity=1 ] (527,178.17) .. controls (527,176.95) and (527.99,175.96) .. (529.21,175.96) .. controls (530.43,175.96) and (531.42,176.95) .. (531.42,178.17) .. controls (531.42,179.39) and (530.43,180.37) .. (529.21,180.37) .. controls (527.99,180.37) and (527,179.39) .. (527,178.17) -- cycle ;
%Straight Lines [id:da003456844616948951]
\draw    (507.42,246.17) -- (560.42,247.17) ;
%Straight Lines [id:da07575528822481248]
\draw    (489.21,213.37) -- (505.21,243.96) ;
%Straight Lines [id:da6044204842375247]
\draw    (491.42,211.17) -- (558.21,249.37) ;
%Straight Lines [id:da8858074686335713]
\draw    (507.42,246.17) -- (572.21,214.37) ;
%Straight Lines [id:da71588752663663]
\draw    (572.21,214.37) -- (558.21,247.17) ;
%Straight Lines [id:da35853282622868]
\draw    (489.21,211.17) -- (531.42,178.17) ;
%Straight Lines [id:da39045746398726044]
\draw    (491.42,211.17) -- (570,212.17) ;
%Straight Lines [id:da4426570669731862]
\draw    (529.21,178.17) -- (507.42,246.17) ;
%Straight Lines [id:da641056668389391]
\draw    (529.21,180.37) -- (558.21,244.96) ;

% Text Node
\draw (143,74) node [anchor=north west][inner sep=0.75pt]   [align=left] {$u_1$};
% Text Node
\draw (197,75) node [anchor=north west][inner sep=0.75pt]   [align=left] {$u_2$};
% Text Node
\draw (157,116) node [anchor=north west][inner sep=0.75pt]   [align=left] {(a) $G_1$};
% Text Node
\draw (495,74) node [anchor=north west][inner sep=0.75pt]   [align=left] {$w_1$};
% Text Node
\draw (549,75) node [anchor=north west][inner sep=0.75pt]   [align=left] {$w_2$};
% Text Node
\draw (511,116) node [anchor=north west][inner sep=0.75pt]   [align=left] {(c) $G_3$};
% Text Node
\draw (523,75) node [anchor=north west][inner sep=0.75pt]   [align=left] {$e$};
% Text Node
\draw (304,74) node [anchor=north west][inner sep=0.75pt]   [align=left] {$v_1$};
% Text Node
\draw (358,75) node [anchor=north west][inner sep=0.75pt]   [align=left] {$v_2$};
% Text Node
\draw (318,116) node [anchor=north west][inner sep=0.75pt]   [align=left] {(b) $G_2$};
% Text Node
\draw (331,36) node [anchor=north west][inner sep=0.75pt]   [align=left] {$v_3$};
% Text Node
\draw (144,228) node [anchor=north west][inner sep=0.75pt]   [align=left] {$u_1$};
% Text Node
\draw (198,229) node [anchor=north west][inner sep=0.75pt]   [align=left] {$u_2$};
% Text Node
\draw (158,270) node [anchor=north west][inner sep=0.75pt]   [align=left] {(A) $\Gamma(I(G_1))$};
% Text Node
\draw (562.42,242.17) node [anchor=north west][inner sep=0.75pt]   [align=left] {$e^{-1}$};
%\draw (562.42,250.17) node [anchor=north west][inner sep=0.75pt]   [align=left] {$e^{-1}$};
% Text Node
\draw (512,270) node [anchor=north west][inner sep=0.75pt]   [align=left] {(C) $\Gamma(I(G_3))$};
% Text Node
\draw (499,249.17) node [anchor=north west][inner sep=0.75pt]   [align=left] {$e$};
% Text Node
\draw (299,228) node [anchor=north west][inner sep=0.75pt]   [align=left] {$v_1$};
% Text Node
\draw (370,229) node [anchor=north west][inner sep=0.75pt]   [align=left] {$v_2$};
% Text Node
\draw (319,270) node [anchor=north west][inner sep=0.75pt]   [align=left] {(B) $\Gamma(I(G_2))$};
% Text Node
\draw (332,190) node [anchor=north west][inner sep=0.75pt]   [align=left] {$v_3$};
% Text Node
\draw (473,192) node [anchor=north west][inner sep=0.75pt]   [align=left] {$w_2$};
% Text Node
\draw (571,194) node [anchor=north west][inner sep=0.75pt]   [align=left] {$ee^{-1}$};
% Text Node
\draw (520,160) node [anchor=north west][inner sep=0.75pt]   [align=left] {$w_1$};
\end{tikzpicture}
\end{center}

\section{Declarations}
%
%{\bf Ethical Approval:} Not applicable.
%\vspace{5mm}
%
{\bf Competing interests:} Non-financial interests that are directly or indirectly related to
the work submitted for publication.

\vspace{5mm}

{\bf Authors' contributions:} Yanhui Wang: Conceptualization,  Writing-original draft, Editing; Xinyi Zhu:
Writing-original draft; Pei Gao: Conceptualization, Writing-original draft. All authors have read and agreed to the published version of the manuscript.
%\vspace{5mm}
%
%{\bf Funding:}  Not applicable.
%\vspace{5mm}
%
%{\bf Availability of data and materials:} Not applicable.
%
\section*{Acknowledgements}

The authors are very grateful to Professor Victoria Gould for her comments during the preparation of this paper.

\end{document}